\documentclass[a4paper,12pt]{article}
\usepackage[latin1]{inputenc}
\usepackage[francais,english]{babel}    
\usepackage{amssymb}
\usepackage{amsmath}
\usepackage{amsthm}
\usepackage{multicol}
\usepackage{graphicx}
\usepackage{color}
\usepackage[T1]{fontenc}
\usepackage{yfonts}
\usepackage{placeins}
\usepackage{array}
\usepackage{here}

\newcommand{\sign}{\begin{flushright}
Thomas Haettel \\
École Normale Supérieure, DMA UMR 8553 CNRS \\
45 rue d'Ulm, 75005 Paris \\
thomas.haettel@ens.fr
\end{flushright}}

\newtheorem{thm}{Théorème}[section]
\newcommand{\bthm}{\begin{thm}}
\newcommand{\ethm}{\end{thm}}

\newtheorem*{thn}{Théorème}
\newcommand{\bthn}{\begin{thn}}
\newcommand{\ethn}{\end{thn}}

\newtheorem{defi}[thm]{Définition}
\newcommand{\bdf}{\begin{defi}}
\newcommand{\edf}{\end{defi}}

\newtheorem{qc}{Question de Cours}
\newcommand{\bqc}{\begin{qc}}
\newcommand{\eqc}{\end{qc}}

\newtheorem{rqc}{R\'{e}ponse}
\newcommand{\brqc}{\begin{rqc}}
\newcommand{\erqc}{\end{rqc}}

\newtheorem{exo}{Exercice}
\newcommand{\bex}{\begin{exo}}
\newcommand{\eex}{\end{exo}}

\newtheorem{sol}{Solution}
\newcommand{\bsol}{\begin{sol}}
\newcommand{\esol}{\end{sol}}

\newtheorem{pro}[thm]{Proposition}
\newcommand{\bpro}{\begin{pro}}
\newcommand{\epro}{\end{pro}}

\newtheorem{cor}[thm]{Corollaire}
\newcommand{\bcor}{\begin{cor}}
\newcommand{\ecor}{\end{cor}}

\newtheorem{lem}[thm]{Lemme}
\newcommand{\blem}{\begin{lem}}
\newcommand{\elem}{\end{lem}}

\newtheorem*{rmq}{Remarque}
\newcommand{\brq}{\begin{rmq} \upshape}
\newcommand{\erq}{\end{rmq}}

\newtheorem*{exe}{Exemple}
\newcommand{\bexe}{\begin{exe} \upshape}
\newcommand{\eexe}{\end{exe}}

\newtheorem*{pre}{Démonstration}
\newcommand{\bp}{\begin{pre} \upshape}
\newcommand{\ep}{\hfill $\square$ \end{pre}}

\newcommand{\beq}{\begin{eqnarray*}}
\newcommand{\eeq}{\end{eqnarray*}}

\newcommand{\bfg}{
\begin{figure}[H]
\begin{center}}
\newcommand{\efg}{
\end{center}
\end{figure}
\FloatBarrier}

\newcolumntype{M}[1]{>{\raggedright}m{#1}}

\newcommand{\df}{\emph}

\newcommand{\R}{\mathbb{R}}
\newcommand{\Q}{\mathbb{Q}}
\newcommand{\N}{\mathbb{N}}
\newcommand{\Z}{\mathbb{Z}}
\newcommand{\C}{\mathbb{C}}

\renewcommand{\P}{\mathbb{P}}
\renewcommand{\SS}{\mathbb{S}}

\newcommand{\bs}{\symbol{92}}
\newcommand{\Ch}{\mathcal{G}}

\newcommand{\ov}{\overline}

\newcommand{\ra}{\rightarrow}

\renewcommand{\geq}{\geqslant}
\renewcommand{\leq}{\leqslant}

\makeatletter
\def\Ddots{\mathinner{\mkern1mu\raise\p@
\vbox{\kern7\p@\hbox{.}}\mkern2mu
\raise4\p@\hbox{.}\mkern2mu\raise7\p@\hbox{.}\mkern1mu}}
\makeatother

\makeatletter
\def\maketitles{%
  \null
  \thispagestyle{empty}%
  \vfill
  \begin{center}\leavevmode
    \normalfont
    {\LARGE \@title\par}%
    \vskip 1.2cm
    {\large \@author\par}%
    \vskip 1.2cm
    {\large \@subtitle\par}%
    \vskip 0.8cm
    {\large \@date\par}%
  \end{center}%
  \vfill
  \null
  \cleardoublepage
  }

\def\date#1{\def\@date{#1}}
\def\author#1{\def\@author{#1}}
\def\title#1{\def\@title{#1}}
\def\subtitle#1{\def\@subtitle{#1}}
\makeatother

\title{L'espace des sous-groupes fermés de $\R \times \Z$}
\author{Thomas Haettel}
\date{}

\begin{document}

\maketitle

\begin{abstract} The space of closed subgroups of a locally compact topological group is endowed with a natural topology, called the Chabauty topology. We completely describe the space of closed sugroups of the group $\R \times \Z$, which is not trivial~: for example, its fundamental group is uncountable.
\end{abstract}

\selectlanguage{francais}

\bigskip

Si $G$ est un groupe topologique localement compact, l'espace $\Ch(G)$ des sous-groupes fermés de $G$ est muni de la \textit{topologie de Chabauty} (voir~\cite{chabauty}). Cette topologie fait de $\Ch(G)$ un espace compact (voir par exemple~\cite{harpe_chabauty} pour une excellente introduction). L'objet de cet article est l'étude de cet espace pour le groupe topologique $G=\R \times \Z$, et nous montrerons que la topologie de $\Ch(\R \times \Z)$ est singulièrement compliquée, et ce malgré la simplicité de $\R \times \Z$.

\medskip

C.~Chabauty a introduit en 1950 cette topologie afin d'étudier les réseaux de $\R^n$. L'espace $\Ch(G)$ dans son ensemble est en général difficile à expliciter : le premier calcul complet non banal est dû à I.~Pourezza et J.~Hubbard en 1979, qui ont montré que pour le groupe topologique $G = \R^2$, l'espace $\Ch(\R^2)$ est homéomorphe à la sphère $\SS^4$ de dimension $4$ (voir~\cite{pourezza}). Bien plus récemment, M.~R.~Bridson, P.~de~la~Harpe et V.~Kleptsyn ont calculé cet espace pour le groupe de Heisenberg de dimension $3$ (voir~\cite{heisenberg}).

\medskip

Dans une première partie, nous rappelons des propriétés élémentaires de la topologie de Chabauty, et des recollements d'espaces topologiques. Dans une deuxième partie, nous commençons par décrire les familles de sous-groupes fermés de $\R \times \Z$, puis nous expliquons la manière dont ces sous-espaces de sous-groupes se recollent.

\medskip

Soit $A \subset \R^2$ l'espace des \og anneaux hawaïens \fg\ (voir la figure~\ref{fig:hawai}), réunion d'une infinité dénombrable de cercles $(A_n)_{n \in \N \bs \{0\}}$ se rencontrant deux à deux exactement en un point et s'accumulant sur ce point. Pour tout entier $k \in \N \bs \{0\}$, considérons une copie $\ov{C_k}$ du cône topologique fermé sur le cercle $\R / \Z$, de sorte que la suite des cônes deux à deux disjoints $(\ov{C_k})_{k \in \N \bs \{0\}}$ s'accumule sur le segment $I=[0,\infty]$ le long des génératrices de ces cônes. Enroulons enfin le bord du cône $\ov{C_k}$ sur l'espace $A$, de sorte que le cercle $A_n$ soit parcouru $0$ fois si $k$ ne divise pas $n$ et $\varphi(\frac{n}{k})$ fois si $k$ divise $n$ (où $\varphi$ désigne la fonction indicatrice d'Euler). Une définition plus précise de l'application de recollement de $\bigcup_{k \in \N \bs \{0\}} \ov{C_k}$ sur $A$ sera donnée dans la seconde partie, et on pourra se reporter à la figure~\ref{fig:espacex} pour une représentation schématique.

\bthn L'espace $\Ch(\R \times \Z)$ est homéomorphe à la réunion des cônes $(\ov{C_k})_{k \in \N \bs \{0\}}$ s'accumulant sur $I$, recollés sur l'espace $A$. \ethn

Enfin, nous nous intéressons au groupe fondamental de l'espace $\Ch(\R \times \Z)$, et nous montrons que malgré le recollement de tous ces disques sur $A$, il n'est pas dénombrable.

Je tiens à remercier chaleureusement Frédéric Paulin pour ses relectures attentives et ses précieux conseils, ainsi que le rapporteur pour ses nombreuses remarques ayant permis de bien améliorer la présentation de cet article.

\section{Préliminaires}

\subsection{Définitions}

Soit $X$ un espace topologique localement compact, et $\mathcal{F}(X)$ l'ensemble des fermés de $X$. On munit $\mathcal{F}(X)$ de la \textit{topologie de Chabauty} (\cite{chabauty}) : les ouverts sont les réunions quelconques d'intersections finies de parties de la forme :
\begin{eqnarray*}
O_K = \{H \in \Ch(G):H \cap K = \emptyset \} \\
O'_U = \{H \in \Ch(G):H \cap U \neq \emptyset \}
\end{eqnarray*}
où $K$ est un compact de $X$ et $U$ un ouvert de $X$.

\medskip

Soit $G$ un groupe topologique localement compact. On note $\Ch(G) \subset \mathcal{F}(G)$ l'ensemble de ses sous-groupes fermés, muni de la topologie induite. Le résultat suivant est classique :

\bpro \label{pro:chg_compact} L'espace topologique $\mathcal{F}(X)$ est compact. De plus, le sous-espace $\Ch(G)$ est un fermé de $\mathcal{F}(G)$, donc est compact. (Voir par exemple~\cite[Chap.~VIII, $\S 5$]{bourbaki}, \cite[Proposition~I.3.1.2, p.~59]{CEG}, \cite[Proposition~1.7, p.~58]{cdp}.) \hfill $\square$ \epro

\bpro \label{pro:ouvert} Soit $f : G \ra H$ un morphisme de groupes topologiques localement compacts, qui soit une application ouverte. Alors l'application $\Ch^*(f) : \Ch(H) \ra \Ch(G)$ définie par $A \mapsto f^{-1}(A) $ est continue. \epro

\bp Soit $K$ un compact de $G$, alors $f(K)$ est un compact de $H$ et $\Ch^*(f)^{-1}(O_K) = O_{f(K)}$ est un ouvert de $\Ch(H)$. Soit $U$ un ouvert de $G$, alors par hypothèse $f(U)$ est un ouvert de $H$ et $\Ch^*(f)^{-1}(O'_U) = O'_{f(U)}$ est un ouvert de $\Ch(H)$. Ainsi l'application $\Ch^*(f)$ est continue. \ep

\bpro \label{pro:surjection} Soit $f : G \ra H$ un morphisme de groupes topologiques localement compacts, qui soit une surjection ouverte. Alors l'application $\Ch^*(f) : \Ch(H) \ra \Ch(G)$ définie par $A \mapsto f^{-1}(A) $ est un homéomorphisme sur son image. \epro

\bp Puisque le morphisme $f$ est surjectif, pour tout sous-groupe fermé $A$ de $H$, nous avons $f(\Ch^*(f)(H))=H$, donc l'application $\Ch^*(f)$ est injective. D'après la proposition précédente, l'application $\Ch^*(f)$ est continue. Enfin, puisque l'espace $\Ch(H)$ est compact et l'espace $\Ch(G)$ séparé, l'application continue injective $\Ch^*(f)$ est un homéomorphisme sur son image. \ep

\bpro \label{pro:injection} Soit $f:G \rightarrow H$ un morphisme de groupes topologiques localement compacts qui est un homéomorphisme sur son image, d'image fermée. Alors l'application $\Ch_*(f) : \Ch(G) \ra \Ch(H)$ définie par $A \mapsto f(A) $ est un homéomorphisme sur son image. \epro

\bp Soit $A$ un sous-groupe fermé de $G$. Alors, puisque $f$ est un homéomorphisme sur son image, $f(A)$ est un sous-groupe fermé de $f(G)$. Or $f(G)$ est lui-même un sous-groupe fermé de $H$, donc finalement $f(A)$ est un sous-groupe fermé de $H$ : ainsi, l'application $\Ch_*(f)$ est bien définie.

Considérons l'application $g = f^{-1}|_{f(G)} : f(G) \ra G$ : c'est un isomorphisme de groupes topologiques localement compacts. Ainsi l'application $\Ch^*(g) : \Ch(G) \ra \Ch(H)$ est un homéomorphisme sur son image, et cette application coïncide avec l'application $\Ch_*(f)$. \ep

Remarquons que l'hypothèse que $f$ réalise un homéomorphisme sur son image est nécessaire. En effet, considérons l'identité $i$ du groupe $G=\R$ muni de la topologie discrète, à valeurs dans $H=\R$ muni de la topologie usuelle. L'identité $i$ est un morphisme de groupes topologiques, bijectif, dont l'image est bien un sous-groupe fermé de $\R$, mais l'application $\Ch_*(i)$ n'est même pas définie : tout sous-groupe (même non fermé) de $H$ est un sous-groupe fermé de $G$.

\bigskip

Dans le cas où le groupe topologique localement compact $G$ est muni d'une distance induisant sa topologie, on peut décrire la convergence des suites de sous-groupes fermés de $G$ pour la topologie de Chabauty.

\bpro Une suite de sous-groupes fermés $(H_n)_{n \in \N}$ converge vers un sous-groupe fermé $H$ dans $\Ch(G)$ si et seulement si $H$ est l'ensemble des valeurs d'adhérence des suites de $(H_n)_{n \in \N}$, c'est-à-dire :
\begin{enumerate}
\item Pour tout $x \in H$, il existe une suite $(x_n)_{n \in \N}$ convergeant vers $x$ telle que, pour tout $n$, nous ayons $x_n \in H_n$.
\item Pour toute partie infinie $P \subset \N$, pour toute suite $(x_n)_{n \in P}$ convergeant vers $x$ telle que $x_n \in H_n$ pour tout $n \in P$, nous avons $x \in H$.
\end{enumerate}
(Voir par exemple~\cite[Lemma~I.3.1.3, p.~60]{CEG}, \cite[Proposition~1.8, p.~60]{cdp}.) \hfill $\square$ \epro

\bpro \label{pro:seq} Si de plus la distance $d$ sur $G$ est propre (i.e. les boules fermées sont compactes), alors l'espace $\Ch(G)$ est métrisable, pour la \df{distance de Hausdorff} pointée (voir~\cite[Definition~5.43, p.~76]{bridson_haefliger}) : si $H$, $H'$ sont des sous-groupes fermés de $G$, on définit $d_{\mbox{\scriptsize Hau}}(H,H')$ comme la borne inférieure des $\varepsilon >0$ tels que :
\beq H \cap B(e,\frac{1}{\varepsilon}) & \subset & V_\varepsilon(H') \\
\mbox{et } H' \cap B(e,\frac{1}{\varepsilon}) & \subset & V_\varepsilon(H),\eeq 
où $V_\varepsilon(H')$ désigne le $\varepsilon$-voisinage ouvert de $H'$ dans $G$.
(Voir par exemple~\cite[Proposition~1.8, p.~60]{cdp}.) \hfill $\square$ \epro

\subsection{Exemples}

Les exemples suivants de calculs d'espaces des sous-groupes fermés sont bien connus, nous les rappelons pour fixer les notations.

\bigskip

Notons $X_\Z$ le sous-espace topologique compact de $\R$ défini par $X_\Z = \{0\} \cup \{\frac{1}{n}, n \in \N \bs \{0\} \}$.

\bpro L'application $\phi_\Z : X_\Z  \rightarrow \Ch(\Z)$ définie par $\frac{1}{n} \mapsto n\Z$ et $0 \mapsto \{0\}$ est un homéomorphisme. \hfill $\square$ \epro

Adoptons la convention suivante désormais : $\frac{1}{0}\Z = \{0\}$ et $\frac{1}{\infty}\Z = \R$. Cette convention est justifiée par la proposition suivante.

\bpro \label{pro:R} L'application $\phi_\R : X_\R=[0,\infty] \rightarrow \Ch(\R)$ définie par $\alpha \mapsto \frac{1}{\alpha}\Z$ est un homéomorphisme. \hfill $\square$ \epro

\subsection{Recollement d'espaces topologiques}

Soient $X,Y$ deux espaces topologiques disjoints, $A \subset X$ une partie fermée et $f:A \rightarrow Y$ une application continue. Soit $\mathcal{R}$ la relation d'équivalence sur $X \sqcup Y$ engendrée par $\{a \sim f(a):a \in A\}$. Alors l'espace topologique quotient $X \cup_f Y = (X \sqcup Y) / \mathcal{R}$ est appelé le \textit{recollement} de $X$ et $Y$ via $f$ le long de $A$.

\bpro \label{pro:recollement_fonction} Soient $X,Y,A,f$ comme ci-dessus et $p:X \sqcup Y \rightarrow X \cup_f Y$ l'application de projection. Soient $Z$ un espace topologique, $\phi:X \rightarrow Z$ et $\psi:Y \rightarrow Z$ continues et compatibles pour $f$, c'est-à-dire qui vérifient $\phi|_A = \psi \circ f$. Alors l'unique application $\phi \cup_f \psi : X \cup_f Y \rightarrow Z$ telle que $(\phi \cup_f \psi) \circ p|_X = \phi$ et $(\phi \cup_f \psi) \circ p|_Y = \psi$ est continue. (Voir par exemple~\cite[Theorem~6.5, p.~129]{dugundji}.) \epro

\bpro \label{pro:recollement_separe} Soient $X,Y,A,f$ comme ci-dessus. Supposons $X$ régulier (i.e. on peut séparer un point et un fermé disjoints par des ouverts) et $Y$ séparé, alors le recollement $X \cup_f Y$ est séparé. \epro

\bp Deux points distincts de $p(Y)$ sont séparés dans $X \cup_f Y$. Soit $x \in X \bs A$ et $a \in A$. Montrons que $p(x)$ et $p(a)$ sont séparés dans $X \cup_f Y$ : soit $F=f^{-1}(f(a))$ fermé de $A$. Or $A$ est fermé dans $X$, donc $F$ est fermé dans $X$. De plus $F$ est inclus dans $A$, et $x \not\in A$, donc $x \not\in F$. Comme $X$ est régulier, on peut séparer $x$ et $F$ par des ouverts de $X$, lesquels se projettent par $p$ en des ouverts disjoints de $X \cup_f Y$. \ep

Soit $X$ un espace topologique et $A$ une partie fermée de $X$. Soit $\mathcal{R}$ la relation d'équivalence engendrée par $\{a \sim a':a,a' \in A\}$. On note $X / \langle A \rangle$ l'espace topologique quotient $X / \mathcal{R}$, appelé l'\textit{écrasement} de $A$ dans $X$. Il est naturellement homéomorphe au recollement $X \cup_f \{pt\}$, où $f:A \rightarrow \{pt\}$ est une application constante. En particulier, d'après la proposition~\ref{pro:recollement_separe}, si $X$ est régulier, alors l'écrasement $X / \langle A \rangle$ est séparé.

\section{Le groupe $\R \times \Z$}

Considérons le groupe topologique localement compact $\R \times \Z$, muni de la topologie produit~: le groupe $\R \times \Z$ est métrisable, pour la distance $\delta((x,n),(x',n')) = \max \{|x-x'|,|n-n'| \}$. Remarquons de plus que cette distance est propre. Notons $i : \R \hookrightarrow \R \times \Z$ le morphisme injectif $x \mapsto (x,0)$ et $\pi:\R \times \Z \rightarrow \Z$ la deuxième projection.

Nous noterons $E(x)$ la partie entière du réel $x$. Adoptons la convention que nous écrirons chaque rationnel $\beta=\frac{a}{b}$ avec $a \in \Z$ et $b \in \N \bs \{0\}$ premiers entre eux. Si $\beta \in \R$, nous noterons $\ov{\beta}$ son image dans $\R/\Z$.

\subsection{Description des sous-groupes fermés de $\R \times \Z$}

Nous allons maintenant décrire tous les sous-groupes fermés du groupe $\R \times \Z$.

\bigskip

Remarquons que le morphisme $i$ est ouvert, donc d'après la proposition~\ref{pro:ouvert}, l'application $\Ch^*(i) : \Ch(\R \times \Z) \ra \Ch(\R)$ qui à $H$ associe $i^{-1}(H)$ est continue. Grâce à l'homéomorphisme $\phi_\R : [0,\infty] \ra \Ch(\R)$, on peut définir l'application continue $\alpha : \phi_\R^{-1} \circ \Ch^*(i) : \Ch(\R \times \Z) \ra [0,\infty]$.

Soit $H$ un sous-groupe fermé de $\R \times \Z$. On notera $\alpha=\alpha(H)$ son image par $\alpha$. De plus, $\pi(H)$ est un sous-groupe de $\Z$ : soit donc $n=n(H)$ l'unique élément de $\N$ tel que $\pi(H)=n\Z$.

\bigskip

Si $n=0$, alors $H=G^I_\alpha$, où nous notons
$$\hspace{0.8cm} G^I_\alpha = \Z(\frac{1}{\alpha},0).$$

Si $n>0$, plusieurs cas sont à distinguer.
\begin{itemize}
\item Si $\alpha=0$, alors $\pi^{-1}(n) \cap H=\{ (\gamma,n) \}$, pour un unique $\gamma=\gamma(H) \in \R$. Alors $H=G^{II}_{\gamma,n}$, où nous notons 
$$G^{II}_{\gamma,n} = \Z(\gamma,n).$$
\item Si $0<\alpha<\infty$, alors $\pi^{-1}(n) \cap H=\{(\frac{\beta + p}{\alpha},n),p \in \Z\}$, pour un unique $\ov{\beta}=\ov{\beta}(H) \in \R / \Z$. Alors $H=G^{III}_{\alpha,\ov{\beta},n}$, où nous notons
$$G^{III}_{\alpha,\ov{\beta},n} = \Z(\frac{1}{\alpha},0)+\Z(\frac{\beta}{\alpha},n).$$
\item Si $\alpha=\infty$, alors $\pi^{-1}(n) \cap H=\R \times \{ n \}$. Alors $H = G^{IV}_n$, où nous notons
$$G^{IV}_n = \R \times n\Z.$$
\end{itemize}

\smallskip

\bpro \label{pro:ChG} L'ensemble $\Ch(\R \times \Z)$ des sous-groupes fermés de $\R \times \Z$ est réunion disjointe des familles
$$\{G^I_\alpha = \Z(\frac{1}{\alpha},0) \,:\, \alpha \in [0,\infty]\}$$
$$\{G^{II}_{\gamma,n} = \Z(\gamma,n) \,:\, \gamma \in \R,n \in \N \bs \{0\} \}$$
$$\{G^{III}_{\alpha,\ov{\beta},n} = \Z(\frac{1}{\alpha},0)+\Z(\frac{\beta}{\alpha},n) \,:\, \alpha \in ]0,\infty[,\ov{\beta} \in \R / \Z,n \in \N \bs \{0\} \}$$
$$\mbox{ et } \{G^{IV}_n = \R \times n\Z \,:\, n \in \N \bs \{0\} \}.$$
De plus, le paramétrage de chacune de ces familles est bijectif. \epro

\bp Ceci découle de la définition et de l'unicité des paramètres $\alpha(H)$, $n(H)$, $\ov{\beta}(H)$ (si $n>0$ et $\alpha \in \;]0,\infty[\;$) et $\gamma(H)$ (si $n>0$ et $\alpha=0$), pour un sous-groupe de $\R \times \Z$ fermé $H$ donné. \ep

Considérons trois sous-espaces de $\Ch(\R \times \Z)$ :
\begin{enumerate}
\item Le sous-espace $\Ch^I$ des sous-groupes fermés de $\R \times \Z$ dont la projection sur $\Z$ est $\{0\}$. Ses éléments sont les sous-groupes $G^I_\alpha$, pour $\alpha \in [0,\infty]$.
\item Le sous-espace $\Ch^{II}$ des sous-groupes fermés de $\R \times \Z$ qui sont cycliques infinis et ont une projection sur $\Z$ différente de $\{0\}$, ainsi que le sous-groupe $\{0\}$. Ses éléments sont les sous-groupes $G^{II}_{\gamma,n}$, pour $\gamma \in \R$ et $n \in \N \bs \{0\}$, et $\{0\}$.
\item Le sous-espace $\Ch^{III}$ des sous-groupes fermés de $\R \times \Z$ isomorphes à $\Z^2$ ou à $\R \times \Z$. Ses élements sont les sous-groupes $G^{III}_{\alpha,\ov{\beta},n}$, pour $\alpha \in ]0,\infty[$, $\ov{\beta} \in \R / \Z$ et $n \in \N \bs \{0\}$, ainsi que les sous-groupes $G^{IV}_n$, pour $n \in \N \bs \{0\}$. Notons $\Ch^{III}_n$ le sous-espace correspondant à une valeur $n \in \N \bs \{0\}$ fixée, c'est-à-dire la réunion des $G^{III}_{\alpha,\ov{\beta},n}$, pour $\alpha \in \,]0,\infty[$ et $\ov{\beta} \in \R / \Z$, et de $G^{IV}_n$.
\end{enumerate}

Nous allons décrire la topologie de chacun de ces sous-espaces $\Ch^I$, $\Ch^{II}$ et $\Ch^{III}$. Puis nous allons décrire comment ces espaces se recollent pour former l'espace $\Ch(\R \times \Z)$.

\subsection{Le sous-espace $\Ch^I$}

\bpro \label{pro:diag} L'application $\psi^I : [0,\infty] \rightarrow \Ch^I$ définie par $\alpha \mapsto \Z(\frac{1}{\alpha},0)$ est un homéomorphisme. \epro

\bp Remarquons que $\psi^I$ est la composée de l'homéomorphisme $\phi_\R : [0,\infty] \ra \Ch(\R)$ et de l'application $\Ch_*(i) : \Ch(\R) \ra \Ch(\R \times \Z)$. Or le morphisme $i$ est un plongement d'image fermée de $\R \times \Z$, donc d'après la proposition~\ref{pro:injection}, l'application $\Ch_*(i)$ est un plongement, d'image $\Ch^I$. \ep

\subsection{Le sous-espace $\Ch^{II}$}

Considérons l'espace topologique des \og anneaux hawaïens \fg
$$A = \bigcup_{n \in \N \bs \{0\} } A_n$$
où $A_n$ désigne le cercle dans la droite complexe $\C$ de centre $\frac{1}{n}$ et de rayon $\frac{1}{n}$ (voir la figure~\ref{fig:hawai}).
\bfg
\includegraphics{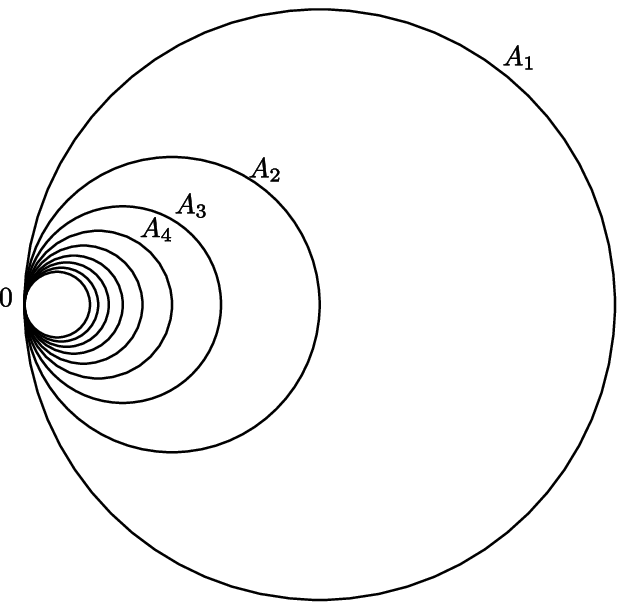}
\caption{L'espace des \og anneaux hawaïens \fg\ A}
\label{fig:hawai}
\efg
Considérons la bijection $\psi^{II}:A  \rightarrow \Ch^{II}$ définie par
\begin{eqnarray*}
\frac{1}{n}(1+e^{2i\theta}) \neq 0 & \mapsto & G^{II}_{n\tan{\theta},n} \mbox{ , où } n \in \N \bs \{0\} \mbox{ et } \theta \in \,]-\frac{\pi}{2},\frac{\pi}{2}[\, \\
0 & \mapsto & \{0\}.
\end{eqnarray*}

\smallskip

L'application $\psi^{II}$ admet pour inverse l'application $(\psi^{II})^{-1}: \Ch^{II}  \rightarrow  A$ définie par
\begin{eqnarray*}
G^{II}_{\gamma,n} & \mapsto & \frac{1}{n}(1+e^{2i\arctan(\frac{\gamma}{n})}) \mbox{ où } n \in \N \bs \{0\} \mbox{ et } \gamma \in \R \\
\{0\} & \mapsto & 0.
\end{eqnarray*}

\bpro \label{pro:hawai} L'application $\psi^{II}$ est un homéomorphisme. \epro

\bp Comme l'espace de départ est compact métrisable et que l'espace d'arrivée est séparé métrisable, il suffit de montrer que l'application $\psi^{II}$ est séquentiellement continue.

Montrons que, pour tout $n \in \N \bs \{0\}$ et $\theta \in \,]-\frac{\pi}{2},\frac{\pi}{2}[$, l'application $\psi^{II}$ est continue en $z=\frac{1}{n}(1+e^{2i\theta})$. Soit $(z_k=\frac{1}{n}(1+e^{2i\theta_k}))_{k \in \N}$ une suite de $A$ convergeant vers $z$, c'est-à-dire que la suite $(\theta_k)_{k \in \N}$ converge vers $\theta$. Montrons que la suite $(\psi^{II}(z_k))_{k \in \N}$ converge vers $\psi^{II}(z)$. Le générateur $(n\tan{\theta},n)$ de $\psi^{II}(z)$ est limite de la suite $(n\tan{\theta_k},n)_{k \in \N}$ d'élements de $(\psi^{II}(z_k))_{k \in \N}$. Réciproquement, supposons qu'une suite $(p_kn\tan{\theta_k},p_kn)_{k \in P}$ d'élements de $(\psi^{II}(z_k))_{k \in P}$ converge vers $(x,m)$, où  $P$ désigne une partie infinie de $\N$ et où $p_k \in \Z$ pour tout $k \in P$. Alors $p_k=p$ est constant à partir d'un certain, donc $(x,m) = p(n\tan{\theta},n) \in \psi^{II}(z)$.

Montrons que l'application $\psi^{II}$ est continue en $0$. Soit $(z_k=\frac{1}{n_k}(1+e^{2i\theta_k}))_{k \in \N}$ une suite de $A$ convergeant vers $0$. Si la suite $(n_k)_{k \in \N}$ tend vers $+\infty$, il est clair que la suite de sous-groupes $(\psi^{II}(z_k))_{k \in \N}$ converge vers $\{0\} = \psi^{II}(0)$. Sinon, quitte à extraire, on peut supposer que la suite $(\theta_k)_{k \in \N}$ tend vers $\frac{\pm \pi}{2}$, et dans ce cas il est également clair que la suite de sous-groupes $(\psi^{II}(z_k))_{k \in \N}$ converge vers $\{0\} = \psi^{II}(0)$. \ep

Définissons de plus, pour tout $n \in \N \bs \{0\}$, l'application $\psi^{II}_n : A \ra \Ch^{II}$ par
\beq
\frac{1}{b}(1+e^{2i\theta}) \neq 0 & \mapsto & G^{II}_{b\tan{\theta},bn} \mbox{ , où } b \in \N \bs \{0\} \mbox{ et } \theta \in \,]-\frac{\pi}{2},\frac{\pi}{2}[\, \\
0 & \mapsto & \{0\}.
\eeq

Il est clair que l'application $\psi^{II}_n$ est un homéomorphisme sur une partie fermée de $\Ch^{II}$.

\subsection{Le sous-espace $\Ch^{III}$}

\label{sec:GIII}

Notons $\ov{C}$ l'écrasement $\{ (\alpha,\ov{\beta}), \alpha \in [0,\infty],\ov{\beta} \in \R / \Z \} / \langle \{\infty\} \times \R / \Z \rangle$, muni de la topologie quotient de la topologie usuelle sur $[0,\infty] \times \R / \Z$ : l'espace $\ov{C}$ est homéomorphe à un cône fermé (on notera indifféremment un point de $[0,\infty] \times \R / \Z$ et son image dans $\ov{C}$). Notons de plus $C$ l'image de $]0,\infty] \times \R / \Z$ dans $\ov{C}$ : l'espace $C$ est homéomorphe à un cône ouvert. Notons également $\partial \ov{C}$ le bord de $C$, c'est-à-dire l'image de $\{0\} \times \R / \Z$ dans $\ov{C}$ : l'espace $\partial \ov{C}$ est homéomorphe au cercle $\R/\Z$.

\bpro \label{pro:cone} Pour tout $n \in \N \bs \{0\}$, l'application $\psi^{III}_n: C \rightarrow \Ch^{III}_n$ induite par
$$ (\alpha,\ov{\beta}) \mapsto \left\{ \begin{array}{ll} G^{III}_{\alpha,\ov{\beta},n} & \mbox{ si } \alpha < \infty \\
                                                                          G ^{IV}_n  & \mbox{ si } \alpha = \infty \end{array} \right. $$
est un homéomorphisme. \epro

\bp L'application $\psi^{III}_n$ est bijective d'après la proposition~\ref{pro:ChG}. Les espaces $C$ et $\Ch^{III}_n$ étant métrisables, montrons que $\psi^{III}_n$ est séquentiellement continue et propre. On en déduira que l'application $\psi^{III}_n$ est continue, bijective et propre, donc un homéomorphisme.

\bigskip

Soit $c=(\alpha,\ov{\beta}) \in C$ tel que $\alpha \neq \infty$. Soit $(c_k)_{k \in \N}=(\alpha_k,\ov{\beta_k})_{k \in \N}$ une suite de $C$ convergeant vers $c$, ce qui signifie que la suite $(\alpha_k)_{k \in \N}$ converge vers $\alpha$ et que la suite $(\ov{\beta_k})_{k \in \N}$ converge vers $\ov{\beta}$. Alors les deux générateurs $(\frac{1}{\alpha},0)$ et $(\frac{\beta}{\alpha},n)$ du groupe $\psi^{III}_n(c)$ sont les limites des suites $(\frac{1}{\alpha_k},0)_{k \in \N}$ et 
$(\frac{\beta_k}{\alpha_k},n)$ de $(\psi^{III}_n(c_k))_{k \in \N}$ respectivement. Réciproquement, soit $P$ une partie infinie de $\N$ et soit $(\frac{p_k+q_k\beta_k}{\alpha_k},q_kn)_{k \in P}$ une suite de $(\psi^{III}_n(c_k))_{k \in P}$ convergeant vers $(x,m)$. Alors $q_k=q$ est constant à partir d'un certain rang, et donc $p_k=p$ également. En conclusion, l'élément $(x,m)=(\frac{p+q\beta}{\alpha},qn)$ appartient à $\psi^{III}_n(c)$. On a donc montré que la suite de sous-groupes 
$(\psi^{III}_n(c_k))_{k \in \N}$ convergeait vers $\psi^{III}_n(c)$.

Soit $c=(\infty,0)$ le sommet du cône $C$. Soit $(c_k)_{k \in \N}=(\alpha_k,\ov{\beta_k})_{k \in \N}$ une suite de $C$ convergeant vers $c$, ce qui signifie que la suite $(\alpha_k)_{k \in \N}$ tend vers $\infty$. Supposons que $\alpha_k \neq \infty$ pour tout $k \in \N$. Choisissons des représentants $\beta_k$ de $\ov{\beta_k}$ bornés. Soit $(x,qn) \in \psi^{III}_n(c)$, où $x \in \R$ et $q \in \Z$. Alors la suite $\left( \frac{E(x\alpha_k)+q\beta_k}{\alpha_k},qn \right)_{k \in \N}$ d'éléments de $(\psi^{III}_n(c_k))_{k \in \N}$ converge vers $(x,qn)$. Puisque $\psi^{III}_n(c_k) \subset \psi^{III}_n(c)$ pour tout $k \in \N$, on en déduit que la suite de sous-groupes $(\psi^{III}_n(c_k))_{k \in \N}$ converge vers $\psi^{III}_n(c)$.

Ainsi l'application $\psi^{III}_n$ est continue.

\bigskip

Soit $(c_k)_{k \in \N}=(\alpha_k,\ov{\beta_k})_{k \in \N}$ une suite sortant de tout compact de $C$. Montrons par l'absurde que la suite $(\psi^{III}_n(c_k))_{k \in \N}$ sort de tout compact de $\Ch^{III}_n$ : supposons quitte à extraire que cette suite converge vers un sous-groupe $H \in \Ch^{III}_n$. Alors, par continuité de l'application $\alpha$, on en déduit que la suite $(\alpha_k=\alpha(\psi^{III}_n(c_k)))_{k \in \N}$ converge vers $\alpha(H) \in \;]0,\infty]$. Or le sous-espace $\{\alpha \geq \frac{\alpha(H)}{2}\}$ du cône $C$ est homéomorphe au cône fermé sur le cercle $\{\alpha = \frac{\alpha(H)}{2}\}$, et est donc compact. Ceci est une contradiction  avec le fait que la suite $(c_k)_{k \in \N}$ sort de tout compact, et qu'en particulier $\alpha_k < \frac{\alpha(H)}{2}$ à partir d'un certain rang. Ainsi la suite $(\psi^{III}_n(c_k))_{k \in \N}$ sort de tout compact de $\Ch^{III}_n$, et l'application $\Ch^{III}_n$ est donc propre. 
\ep

\subsection{Domination de l'espace $\Ch(\R \times \Z)$}

Les motivations qui nous ont mené à cette construction proviennent des éclatements en géométrie algébrique, même si dans notre cas il s'agit de \og demi - éclatements\fg , réalisés sur un ensemble non algébrique (les points rationnels de $\R \times \{0\}$). La construction s'inspire également beaucoup des méthodes d'éclatement \og à la Denjoy \fg (voir la partie~\ref{sec:denjoy}).

\bigskip

Considérons l'espace topologique produit $Y^0 = [0,\infty] \times \R \times \P_1(\R)^\Q$ : en tant que produit dénombrable d'espaces compacts métrisables et de l'espace localement compact métrisable $\R$, l'espace $Y^0$ est lui-même localement compact métrisable. Nous noterons sous la forme $(\alpha,\beta,(D_r)_{r \in \Q})$ les éléments de $Y^0$, où pour chaque rationnel $r \in \Q$, $D_r$ est une droite affine du plan euclidien orienté $\R^2$ passant par le point $(0,r)$.

\bigskip

Par l'adjectif \og vertical \fg , nous qualifierons toute droite de $\R^2$ parallèle à la droite d'équation $x=0$, et par l'adjectif \og horizontal \fg , nous qualifierons toute droite de $\R^2$ parallèle à la droite d'équation $y=0$. De plus, si $(\alpha,\beta,(D_r)_{r \in \Q}) \in Y^0$, on notera pour tout $r \in \Q$ par $\theta_r \in \;]-\frac{\pi}{2},\frac{\pi}{2}]$ l'angle orienté entre la droite $D_r$ et la droite horizontale : si la droite $D_r$ est verticale, on pose $\theta_r = \frac{\pi}{2}$. On notera de même $(\theta'_r)_{r \in \Q}$ (resp. $(\theta''_r)_{r \in \Q}$) les angles correspondants pour un élément $(\alpha',\beta',(D'_r)_{r \in \Q})$ (resp. $(\alpha'',\beta'',(D''_r)_{r \in \Q})$) de $Y^0$.

\bigskip

Considérons le sous-espace $Y$ de $Y^0$ constitué des éléments $(\alpha,\beta,(D_r)_{r \in \Q})$ de $Y^0$ dont toutes les droites $(D_r)_{r \in \Q}$ sont concourantes en le point $(\alpha,\beta)$ du plan euclidien si $\alpha \neq \infty$, et dont toutes les droites $(D_r)_{r \in \Q}$ sont horizontales si $\alpha=\infty$. Le sous-espace $Y$ est fermé dans $Y^0$, donc est localement compact et métrisable.

\bigskip

Considérons l'espace quotient $Z$ de $Y$ par la relation d'équivalence $y=(\alpha,\beta,(D_r)_{r \in \Q}) \sim y'=(\alpha',\beta',(D'_r)_{r \in \Q})$ si l'on se trouve dans l'un des cas suivants :
\begin{enumerate}
\item $\alpha=\alpha'=0$ et les droites $D_r$ et $D'_r$ sont verticales pour tout $r \in \Q$.
\item $\alpha=\alpha'=0$, $\beta=\frac{a}{b} \in \Q$, $\beta'=\frac{a'}{b'} \in \Q$, $b=b'$ et les droites $D_\beta$ et $D'_{\beta'}$ sont parallèles.
\item $\alpha=\alpha' \in \;]0,\infty[$ et il existe $p \in \Z$ tel que, pour tout $r \in \Q$, les droites $D_r$ et $D'_{r+p}$ soient parallèles. Remarquons qu'alors $\beta' = \beta + p$.
\item $\alpha=\alpha'=\infty$.
\end{enumerate}

C'est bien une relation d'équivalence, et nous noterons dorénavant $[\alpha,\beta,(D_r)_{r \in \Q}]$ la classe d'équivalence de $(\alpha,\beta,(D_r)_{r \in \Q})$.

\bpro L'espace $Z$ est métrisable et compact. \epro

\bp 
Montrons que l'espace $Z$ est séparé : considérons les trois sous-espaces saturés $Y_0=\{(\alpha,\beta,(D_r)_{r \in \Q}) \in Y : \alpha=0\}$, $Y_{]0,\infty[}=\{(\alpha,\beta,(D_r)_{r \in \Q}) \in Y : \alpha \in \,]0,\infty[\}$ et $Y_\infty=\{(\alpha,\beta,(D_r)_{r \in \Q}) \in Y : \alpha=\infty\}$ de $Y$. Deux points de $Y$ n'appartenant pas à un même de ces trois sous-ensembles peuvent être séparés par des ouverts saturés disjoints, à l'aide du paramètre $\alpha$. Considérons donc le cas où ils appartiennent à un même de ces trois sous-ensembles.

\begin{enumerate}
\item Soient $y=(0,\beta,(D_r)_{r \in \Q}), y'=(0,\beta',(D'_r)_{r \in \Q}) \in Y_0$ non équivalents.
\begin{itemize}
\item Supposons que les droites $D'_r$ soient verticales pour tout $r \in \Q$. Alors $\beta=\frac{a}{b} \in \Q$ et la droite $D_{\beta}$ n'est pas verticale : elle fait un angle $\theta_{\beta} \in \;]-\frac{\pi}{2},\frac{\pi}{2}[$ avec l'horizontale. Soit $\varepsilon = \frac{1}{2}(\frac{\pi}{2} - |\theta_{\beta}|) > 0$. Soient \beq U&=&\{(\alpha'',\beta'',(D''_r)_{r \in \Q}) \in Y : \exists p \in \Z, |\theta''_{\frac{p}{b}}-\theta_{\beta}| < \varepsilon\} \\
\mbox{et } U'&=&\{(\alpha'',\beta'',(D''_r)_{r \in \Q}) \in Y : \alpha'' < \varepsilon \mbox{ et } \forall p \in \Z, |\theta''_{\frac{p}{b}}-\theta_{\beta}| > \varepsilon\}. \eeq
L'espace $U'$ est ouvert car la première condition $\alpha'' < \varepsilon$ implique que la seconde condition $\forall p \in \Z, |\theta''_{\frac{p}{b}}-\theta_{\beta}| > \varepsilon$ ne porte que sur un nombre fini d'entiers $p \in \Z$.
Alors $U$ contient $y$, $U'$ contient $y'$, et $U$ et $U'$ sont des ouverts saturés disjoints.
\item Supposons que les droites $(D'_r)_{r \in \Q}$ ne soient pas toutes verticales, ni les droites $(D_r)_{r \in \Q}$. Ceci signifie que $\beta=\frac{a}{b} \in \Q$, $\beta'=\frac{a'}{b'} \in \Q$ et que les angles $\theta_\beta$ et $\theta'_{\beta'}$ appartiennent à $]-\frac{\pi}{2},\frac{\pi}{2}[$. Puisque $y$ et $y'$ ne sont pas équivalents, cela signifie que $b \neq b'$ ou $\theta_\beta \neq \theta'_{\beta'}$. Soient $\varepsilon>0$ et $\eta>0$, définissons \beq U&=&\{(\alpha'',\beta'',(D''_r)_{r \in \Q}) \in Y : \exists p \in \Z, |\theta''_{\frac{p}{b}}-\theta_{\beta}| < \varepsilon \mbox{ et } \alpha'' < \eta\} \\
\mbox{et }U'&=&\{(\alpha'',\beta'',(D''_r)_{r \in \Q}) \in Y : \exists p' \in \Z, |\theta''_{\frac{p'}{b'}}-\theta'_{\beta'}| < \varepsilon \mbox{ et } \alpha'' < \eta\}. \eeq
Pour $\varepsilon>0$ et $\eta>0$ suffisamment petits, les ensembles $\{(\alpha'',\beta'',(D''_r)_{r \in \Q}) \in Y : |\theta''_{\frac{p}{b}}-\theta_{\beta}| < \varepsilon \mbox{ et } \alpha'' < \eta\}$, pour $p \in \Z$, sont disjoints. Ainsi l'entier $p \in \Z$ intervenant dans la définition d'un élément de $U$ est uniquement déterminé (et de même pour $U'$). 

Supposons $b \neq b'$. Si un point $(\alpha'',\beta'',(D''_r)_{r \in \Q})$ appartient à $U \cap U'$, alors il existe $p,p' \in \Z$ tels que $|\theta''_{\frac{p}{b}}-\theta_{\beta}| < \varepsilon$ et $|\theta''_{\frac{p'}{b'}}-\theta'_{\beta'}| < \varepsilon$. Par ailleurs
$$ \tan \theta''_{\frac{p}{b}} = \frac{\beta''-\frac{p}{b}}{\alpha''} \mbox{ et } \tan \theta''_{\frac{p'}{b'}} = \frac{\beta''-\frac{p'}{b'}}{\alpha''}.$$
Ainsi $\alpha'' = \frac{\frac{p'}{b'}-\frac{p}{b}}{\tan \theta''_{\frac{p}{b}} - \tan \theta''_{\frac{p'}{b'}}}$. Or $|\tan \theta''_{\frac{p}{b}} - \tan \theta''_{\frac{p'}{b'}}|$ est majoré par une constante $c$ dépendant de $\theta_{\beta}$, $\theta'_{\beta'}$ et $\varepsilon$. Donc il suffit de choisir $\eta < \frac{1}{bb'c}$ pour que les ensembles $U$ et $U'$ soient disjoints.

Supposons maintenant $b=b'$ et $\theta_\beta \neq \theta'_{\beta'}$. Si un point $(\alpha'',\beta'',(D''_r)_{r \in \Q})$ appartient à $U \cap U'$, alors il existe $p,p' \in \Z$ tels que $|\theta''_{\frac{p}{b}}-\theta_{\beta}| < \varepsilon$ et $|\theta''_{\frac{p'}{b}}-\theta'_{\beta'}| < \varepsilon$. Pourvu que $\varepsilon$ soit suffisamment petit, on doit avoir $p \neq p'$. Or $\alpha'' = \frac{\frac{p'}{b}-\frac{p}{b}}{\tan \theta''_{\frac{p}{b}} - \tan \theta''_{\frac{p'}{b}}}$, et $|\tan \theta''_{\frac{p}{b}} - \tan \theta''_{\frac{p'}{b}}|$ est majoré par une constante $c$ dépendant de $\theta_{\beta}$, $\theta'_{\beta'}$ et $\varepsilon$, il suffit de choisir $\eta < \frac{1}{b^2c}$ pour que les ensembles $U$ et $U'$ soient disjoints.

Dans chacun des deux cas, on a ainsi construit des voisinages ouverts saturés disjoints de $y$ et $y'$.
\end{itemize}
\item Soient $y=(\alpha,\beta,(D_r)_{r \in \Q}), y'=(\alpha',\beta',(D'_r)_{r \in \Q}) \in Y_{]0,\infty[}$ non équivalents. En particulier $\alpha \neq \alpha'$ ou $\beta \neq \beta' \mod \Z$. Dans chacun de ces cas, il est aisé de construire des voisinages saturés disjoints de $y$ et de $y'$.
\item Tous les points de $Y_\infty$ sont équivalents.
\end{enumerate}

Considérons le sous-espace compact $Y'=\{(\alpha,\beta,(D_r)_{r \in \Q}) \in Y \;:\; \beta \in [0,1] \}$ de $Y$. Ainsi l'espace séparé $Z$, image continue du compact $Y'$ par la projection canonique $pr:Y \ra Z$, est compact.

\bigskip

Considérons les trois sous-espaces $Z_0$, $Z_{]0,\infty[}$ et $Z_\infty$ de $Z$, images de $Y_0$, $Y_{]0,\infty[}$ et $Y_\infty$ respectivement par la projection $pr$.

L'ouvert $Z_{]0,\infty[}$, homéomorphe au produit $\SS^1 \times ]0,\infty[$, est à base dénombrable : soit $(U_n)_{n \in \N}$ une base dénombrable d'ouverts de $Z_{]0,\infty[}$. Considérons les ouverts de $Z$
$$ V_{n,b,\theta} = pr(\{(\alpha,\beta,(D_r)_{r \in \Q}) \in Y \;:\; \exists p \in \Z, |\theta_{\frac{p}{b}}-\theta| < \frac{1}{n} \mbox{ et } \alpha < \frac{1}{n} \}),$$
pour $n \in \N \bs \{0\}$, $b \in \N \bs \{0\}$ et $\theta \in \Q\, \cap \,]\frac{-\pi}{2},\frac{\pi}{2}[$, ainsi que
$$ V'_n = pr(\{(\alpha,\beta,(D_r)_{r \in \Q}) \in Y \;:\; \forall p \in \Z, \frac{\pi}{2} - |\theta_{\frac{p}{n}}| < \frac{1}{n} \mbox{ et } \alpha < \frac{1}{n} \}),$$
pour $n \in \N \bs \{0\}$ et enfin
$$ W_n = pr(\{(\alpha,\beta,(D_r)_{r \in \Q}) \in Y \;:\; \alpha > n \}),$$
pour $n \in \N \bs \{0\}$.

Alors les familles $(U_n)$, $(V_{n,b,\theta})$, $(V'_n)$ et $(W'_n)$ forment une base dénombrable d'ouverts de $Z$.  Ainsi l'espace topologique $Z$ est à base dénombrable, or il est compact, donc il est métrisable.
\ep

Nous allons voir comment les espaces $A$ et $C$ se plongent dans l'espace $Z$.

\bpro L'application $i_A : A  \ra  Z$ définie par
\beq
\frac{1}{b}(1+e^{2i\theta}) \neq 0 & \mapsto & [0,\frac{1}{b},(D_r)_{r \in \Q}], \mbox{ où pour tout $r \in \Q$ non congru à $\frac{1}{b}$ modulo $1$,} \\
& & \mbox{$D_r$ est verticale, et où $\theta_{\frac{1}{b}}=\theta$}\\
0 & \mapsto &  [0,0,(D_r)_{r \in \Q}], \mbox{ où pour tout $r \in \Q$, $D_r$ est verticale} \eeq

est un plongement.
\epro

\bp Cette application est continue en $0$ : si $(z_n)_{n \in \N}=(\frac{1}{b_n}(1+e^{2i\theta_n}))_{n \in \N}$ converge vers $0$, quitte à extraire, soit $(b_n)_{n \in \N}$ tend vers $+\infty$, soit $(\theta_n)_{n \in \N}$ tend vers $\pm \frac{\pi}{2}$. Notons la suite des images $i_A(z_n) = [0,\frac{1}{b_n},(D_{r,n})_{r \in \Q}]$. Dans le cas où la suite $(b_n)_{n \in \N}$ tend vers $+\infty$, alors pour tout $r \in \Q$ la suite $(\theta_{r,n})_{n \in \N}$ converge vers $\pm \frac{\pi}{2}$, donc la suite $(i_A(z_n))_{n \in \N}$ converge vers le point $i_A(0)$. Dans le cas où la suite $(\theta_n)_{n \in \N}$ tend vers $\pm \frac{\pi}{2}$, alors la suite $(i_A(z_n))_{n \in \N}$ converge également vers le point $i_A(0)$. L'application $i_A$ est continue et injective. Or l'espace $A$ est compact, et l'espace $Z$ séparé, donc c'est un plongement. \ep

\bpro L'application $i_C : C  \ra  Z$ définie par
\beq 
(\alpha,\ov{\beta}) \;\;\mbox{(où $\alpha \neq \infty$)} & \mapsto &  [\alpha,\beta,(D_r)_{r \in \Q}], \mbox{ où pour tout $r \in \Q$, la droite} \\
& & \mbox{ $D_r$ passe par les points $(0,r)$ et $(\alpha,\beta)$} \\
(\infty,0) & \mapsto &  [\infty,0,(D_r)_{r \in \Q}], \mbox{ où pour tout $r \in \Q$, la droite $D_r$ est horizontale} \eeq
est un plongement.
\epro

\bp Cette application est continue et injective. Par ailleurs, son inverse $i_C^{-1} : i_C(C) \ra C$ est également continue. \ep

On peut remarquer que l'ensemble $Z$ est réunion disjointe des ensembles $i_A(A)$ et $i_C(C)$. 

\bigskip

Fixons un entier $n \in \N \bs \{0\}$, et considérons l'application
\beq \phi_n : Z & \ra & \Ch(G) \\
i_A(a) & \mapsto & \psi^{II}_n(a) \\
i_C(c) & \mapsto & \psi^{III}_n(c). \eeq

\bpro L'application $\phi_n$ est continue. \epro

\bp D'après les propositions~\ref{pro:hawai} et \ref{pro:cone}, on sait que l'application $\phi_n$, en restriction à $i_A(A)$ et à $i_C(C)$, est continue. Puisque $i_C(C)$ est ouvert dans $Z$ (comme complémentaire du compact $i_A(A)$), il suffit de montrer que si une suite $(z_k=[\alpha_k,\beta_k,(D_{k,r})_{r \in \Q}])_{k \in \N}$ de $i_C(C)$ converge vers $z=[\alpha,\beta,(D_r)_{r \in \Q}] \in i_A(A)$, alors la suite $(\psi^{III}_n \circ i_C^{-1}(z_k))_{k \in \N}$ converge vers $\psi^{II}_n \circ i_A^{-1}(z)$. Choisissons des représentants de $z_k$ tels que la suite $(\beta_k)_{k \in \N}$ converge vers $\beta$.

\bigskip

\begin{itemize}
\item Si $\psi^{II} \circ i_A^{-1}(z) = \{0\}$, cela signifie que pour tout $r \in \Q$, la direction de la droite $D_{k,r}$ tend vers la verticale. Remarquons que l'angle $\theta_{k,r}$ entre la droite $D_{k,r}$ et l'horizontale vaut $\theta_{k,r} = \arctan \left( \frac{\beta_k-r}{\alpha_k} \right)$. Ainsi, pour tout $r \in \Q$, la suite $\left(\left| \frac{\beta_k-r}{\alpha_k} \right|\right)_{k \in \N}$ tend vers $+\infty$. Soit $P$ une partie infinie de $\N$ et des entiers $p_k,q_k \in \Z$ tels que la suite $(\frac{p_k+q_k\beta_k}{\alpha_k},q_kn)_{k \in P}$ d'éléments de $(\psi^{III}_n \circ i_C^{-1}(z_k))_{k \in P}$ converge vers $(x,m)$. Alors la suite $q_k$ est constante égale à $q \in \Z$ à partir d'un certain rang.

Supposons $q \neq 0$, alors la suite $\left(\left| \frac{\beta_k-\frac{-p_k}{q}}{\alpha_k} \right|\right)_{k \in P}$ converge vers $\frac{x}{q}$. Or cette suite est minorée par la suite $\left(\left| \frac{\beta_k-\frac{p'_k}{q}}{\alpha_k} \right|\right)_{k \in P}$, où $p'_k$ est l'entier tel que $\left|\beta_k-\frac{p'_k}{q} \right|$ soit minimal : puisque la suite $(\beta_k)_{k \in \N}$ converge, la suite $(p'_k)_{k \in P}$ est bornée. Or, pour tout $r \in \Q$, la suite $\left(\left| \frac{\beta_k-r}{\alpha_k} \right|\right)_{k \in \N}$ tend vers $+\infty$, donc la suite $\left(\left| \frac{\beta_k-\frac{p'_k}{q}}{\alpha_k} \right|\right)_{k \in P}$ tend vers $+\infty$ : c'est une contradiction, puisque cette suite est bornée. Ainsi $q=0$ et $m=0$, puis comme $(\alpha_k)_{k \in \N}$ tend vers $0$, nous avons $x=0$.

En conclusion, la suite $(\psi^{III}_n \circ i_C^{-1}(z_k))_{k \in \N}$ converge vers $\{0\} = \psi^{II}_n \circ i_A^{-1}(z)$.

\item Si $\psi^{II} \circ i_A^{-1}(z) \neq \{0\}$, alors les droites $(D_r)_{r \in \Q}$ ne sont pas toutes verticales : ainsi $\beta = \frac{a}{b} \in \Q$, et c'est la droite $D_\beta$ qui n'est pas verticale. De plus, pour tout $r \in \Q$ différent de $\beta$, la direction de la droite $D_{k,r}$ tend vers la verticale, et l'angle $\theta_{\beta,r} = \arctan \left( \frac{\beta_k-\beta}{\alpha_k} \right)$ converge vers $\theta \in \;]-\frac{\pi}{2},\frac{\pi}{2}[$. Alors la suite $(\frac{-a+b\beta_k}{\alpha_k},bn)_{k \in \N}$ converge vers $(b\tan{\theta},bn)$, qui est un générateur de $\psi^{II}_n \circ i_A^{-1}(z) = G^{II}_{b\tan{\theta},bn}$.

Réciproquement, soit $P$ une partie infinie de $\N$ et des entiers $p_k,q_k \in \Z$ tels que la suite $(\frac{p_k+q_k\beta_k}{\alpha_k},q_kn)_{k \in P}$ d'éléments de $(\psi^{III}_n \circ i_C^{-1}(z_k))_{k \in P}$ converge vers $(x,m)$. Alors la suite $q_k$ est constante égale à $q \in \Z$ à partir d'un certain rang, et la suite $(p_k+q\beta_k)_{k \in P}$ converge vers $0$. Puisque la suite $(\beta_k)_{k \in P}$ converge vers $\beta$, la suite $(p_k)_{k \in P}$ doit être constante à partir d'un certain rang, égale à $p \in \Z$ tel que $\beta = -\frac{p}{q}$. Ainsi il existe $\ell \in \Z$ tel que $q=\ell b$ et $-p = \ell a$. Ainsi, on en conclut que $(x,m) = \ell (b\tan{\theta},bn) \in G^{II}_{b\tan{\theta},bn}$.
\end{itemize}

\bigskip

En conclusion, la suite $(\psi^{III}_n \circ i_C^{-1}(z_k))_{k \in \N}$ converge vers $G^{II}_{b\tan{\theta},bn} = \psi^{II}_n \circ i_A^{-1}(z)$. \ep

Considérons de plus l'application continue
\beq \phi_0 : Z & \ra & \Ch(\R \times \Z) \\
z=[\alpha,\beta,(D_r)_{r \in \Q}] & \mapsto & G^I_{\alpha}. \eeq

Considérons l'espace $X = \Ch(\Z) \times Z$, muni de la topologie produit : c'est un espace métrisable compact. Considérons alors l'application
\beq \phi : X & \ra & \Ch(\R \times \Z) \\
(n\Z , z) & \mapsto & \phi_n(z). \eeq

\bpro L'application $\phi$ est continue et surjective. \epro

\bp La surjectivité est claire. La continuité de chacun des $\phi_n$ assure que l'application $\phi$ est continue sur $(\Ch(\Z) \bs \{\{0\}\}) \times Z$. Montrons la continuité sur $\{\{0\}\} \times Z$.

Soit $(\{0\},z) \in X$, et soit $(n_k\Z,z_k)_{k \in \N}$ une suite de $X$ convergeant vers $(\{0\},z)$. Alors $(\alpha(z_k))_{k \in \N}$ converge vers $\alpha(z)$. Si $n_k=0$ à partir d'un certain rang, alors la continuité de l'application $\phi_0$ assure que la suite de sous-groupes $(\phi_0(z_k))_{k \in \N}$ converge vers $\phi_0(z)$. Quitte à extraire, on peut donc supposer que $n_k>0$ pour tout entier $k \in \N$.
\begin{itemize}
\item Si $\alpha(z)>0$, alors la suite $(n_k)_{k \in \N}$ tend vers $+\infty$. Ainsi la limite de la suite de sous-groupes $(\phi_{n_k}(z_k))_{k \in \N}$ est la même que la suite de sous-groupes $(\Z(\frac{1}{\alpha_k},0))_{k \in \N}$, c'est-à-dire $\Z(\frac{1}{\alpha},0) = \phi_0(z)=\phi(\{0\},z)$.
\item Si $\alpha(z)=0$ et $\alpha(z_k)>0$ à partir d'un certain rang, alors la suite $(z_k)_{k \in \N}$ de $i_C(C)$ converge vers $i_A(0) \in i_A(A)$ dans $Z$, donc pour tout rationnel $r \in \Q$ la suite $(\theta_{k,r})_{k \in \N}$ converge vers $\pm\frac{\pi}{2}$. Soit $P$ une partie infinie de $\N$ et $(p_k)_{k \in P}$, $(q_k)_{k \in P}$ des entiers tels que la suite $(\frac{p_k+q_k\beta_k}{\alpha_k},q_kn_k)_{k \in P}$ d'éléments de $\phi_{n_k}(z_k)$ converge vers $(x,m)$. Ainsi $q_kn_k=m$ à partir d'un certain rang, donc quitte à extraire on peut supposer $q_k=q$ et $n_k=n$ pour tout entier $k \in \N$. On peut de plus supposer que la suite $(\beta_k)_{k \in P}$ converge vers $\beta \in \R$. Ceci implique que la suite $p_k$ doit être constante égale à $p$ à partir d'un certain rang. Si $(p,q) \neq (0,0)$ et $\beta \not\in \Q$, alors $(\frac{p+q\beta_k}{\alpha_k})_{k \in P}$ tend vers $\infty$. Si $(p,q) \neq (0,0)$ et $\beta \in \Q$, alors $\left(\frac{p+q\beta_k}{\alpha_k}=q\tan \theta_{k,-\frac{p}{q}}\right)_{k \in P}$ tend également vers $\infty$. Or cette suite converge vers $x$, donc $(p,q)=(0,0)$ et ainsi $(x,m)=(0,0)$. On a donc montré que la suite de sous-groupes $(\phi_{n_k}(z_k))_{k \in \N}$ convergeait vers $\{0\} = \phi_0(z)$.
\item Si $\alpha(z)=0$ et $\alpha(z_k)=0$ à partir d'un certain rang. Si la suite $(n_k)_{k \in \N}$ tend vers $+\infty$, alors la suite de sous-groupes $(\phi_{n_k}(z_k)=\Z(\gamma_k,n_k))_{k \in \N}$ converge vers $\{0\}=\phi_0(z)$. Sinon on peut supposer que la suite $(n_k)_{k \in \N}$ est constante égale à $n$ à partir d'un certain rang. Dans ce cas, d'après la continuité de l'application $\phi_n$, la suite $(\phi_n(z_k))_{k \in \N}$ converge vers $\phi_n(i_A(0))= \{0\} = \phi_0(z)$.
\end{itemize} \ep

On a ainsi construit un espace compact qui domine l'espace $\Ch(\R \times \Z)$. L'espace $\Ch(\R \times \Z)$ est donc homéomorphe au quotient $\widetilde{X}$ de l'espace $X$ par la relation d'équivalence $(n\Z , z) \sim (n'\Z , z')$ si $\phi(n\Z , z)=\phi(n'\Z , z')$.

\bigskip

Décrivons plus précisément cette relation d'équivalence : deux éléments $(n\Z , z)$ et $(n'\Z , z')$ sont identifiés si l'on se trouve dans l'un de ces trois cas, où l'on note $z=[\alpha,\beta,(D_r)_{r \in \Q}]$ et $z'=[\alpha',\beta',(D'_r)_{r \in \Q}]$.

\begin{enumerate}
\item $\alpha=\alpha'=0$ et les droites $(D_r)_{r \in \Q} = (D'_r)_{r \in \Q}$ sont toutes verticales.
\item $\alpha=\alpha'=0$, $\beta=\frac{a}{b} \in \Q$, $\beta'=\frac{a'}{b'} \in \Q$, $bn=b'n'$ et $\frac{\tan(\theta_\beta)}{n} = \frac{\tan(\theta'_{\beta'})}{n'}$.
\item $\alpha=\alpha'$ et $n=n'=0$.
\end{enumerate}

Afin d'avoir une idée plus précise de la topologie de l'espace $\Ch(\R \times \Z)$, nous allons la décrire avec des recollements.

\subsection{Description avec des recollements}

\label{sec:denjoy}

Considérons le cercle $\R/\Z$, et procédons \og à la Denjoy \fg\ (voir par exemple~\cite[\S~2, p.~403]{katok}) en remplaçant chaque rationnel $\ov{\frac{a}{b}}$ de $\R/\Z$ par un segment $I_{\ov{\frac{a}{b}}}$ : l'espace $R$ ainsi obtenu est encore homéomorphe à un cercle. Une manière de le construire explicitement est de poser $R = \{(\ov{\beta},\theta) \in \R/\Z \times [-\frac{\pi}{2},\frac{\pi}{2}] \;:\; \mbox{si } \ov{\beta} \not\in \Q/\Z \mbox{ alors } \theta = \frac{\pi}{2}\}$, muni de la topologie induite de l'ordre lexicographique sur $\R \times [-\frac{\pi}{2},\frac{\pi}{2}]$ : $(\beta,\theta) \leq (\beta',\theta')$ si $\beta < \beta'$ ou $\beta=\beta'$ et $\theta \leq \theta'$. Fixons un homéomorphisme $\xi=(\xi_1,\xi_2) : \R/\Z \ra R$.

\bigskip

Considérons $\ov{C}$, le cône fermé  sur $\R / \Z$ défini dans la partie~\ref{sec:GIII}, et considérons l'application continue
\beq g : \partial \ov{C} & \ra & A \\
(0,\ov{\beta}) & \mapsto & \left\{ \begin{array}{ll} \frac{1}{b}(1+e^{2i\xi_2(\ov{\beta})}) & \mbox{ si } \xi_1(\ov{\beta}) = \ov{\frac{a}{b}} \in \Q/\Z \\
0 & \mbox{ sinon.} \end{array} \right. \eeq

\bpro L'espace $Z$ est homéomorphe au recollement $\ov{C} \cup_g A$ des espaces $\ov{C}$ et $A$ le long de $g$. \epro

\bp
\
\begin{figure}[H]
\begin{minipage}{0.45\textwidth}
\par Considérons, pour chaque rationnel $\ov{\beta} = \ov{\frac{a}{b}} \in \Q/\Z$, le demi-disque ouvert $\Gamma_{\ov{\beta}}$ inclus dans $C \subset \ov{C}$, de diamètre le segment $\xi^{-1}(\{\ov{\beta}\} \times [-\frac{\pi}{2},\frac{+\pi}{2}])$ : voir la figure~\ref{fig:gamma}. 
\end{minipage}
\hfill
\begin{minipage}{0.45\textwidth}
\centering
\input{gamma.pstex_t}
\caption{Les demi-disques $(\Gamma_{\ov{\beta}})_{\ov{\beta} \in \Q/\Z}$}
\label{fig:gamma}
\end{minipage}
\end{figure}
Notons $Y_0=\{(\alpha,\beta,(D_r)_{r \in \Q}) \in Y \;:\; \alpha=0\}$ et $Y_{>0}=\{(\alpha,\beta,(D_r)_{r \in \Q}) \in Y \;:\; \alpha>0\}$ les deux sous-espaces de $Y$, et $Z_0$ et $Z_{>0}$ leurs images dans $Z$. Identifions $Z_{>0}$ avec le cône ouvert $C$ par l'application $[\alpha,\beta,(D_r)_{r \in \Q}] \mapsto (\alpha,\beta)$.
\
\begin{figure}[H]
\begin{minipage}{0.45\textwidth}
Considérons, pour chaque rationnel $\ov{\beta} = \ov{\frac{a}{b}} \in \Q/\Z$, le disque ouvert $\Gamma'_{\ov{\beta}}$ dans $Z_{>0}$ tangent au cercle $\{\alpha=0\}$ en le point $(0,\beta)$, et de diamètre $\frac{1}{10b^2}$ suffisament petit pour que les disques $(\Gamma'_{\ov{\beta}})_{\ov{\beta} \in \Q/\Z}$ soient deux à deux disjoints : voir la figure~\ref{fig:gamma'}.
\end{minipage}
\hfill
\begin{minipage}{0.45\textwidth}
\centering
\input{gammap.pstex_t}
\caption{Les disques $(\Gamma'_{\ov{\beta}})_{\ov{\beta} \in \Q/\Z}$}
\label{fig:gamma'}
\end{minipage}
\end{figure}

Il est clair qu'il existe un homéomorphisme entre l'espace $C$ privé des demi-disques $(\Gamma_{\ov{\beta}})_{\ov{\beta} \in \Q/\Z}$ et l'espace $Z_{>0}$ privé des disques $(\Gamma'_{\ov{\beta}})_{\ov{\beta} \in \Q/\Z}$, qui envoie le demi-cercle $\partial \Gamma_{\ov{\beta}}$ au bord du demi-disque $\Gamma_{\ov{\beta}}$ sur le cercle $\partial \Gamma'_{\ov{\beta}}$ au bord du disque $\Gamma'_{\ov{\beta}}$. Prolongeons cet homéomorphisme en envoyant le demi-disque $\Gamma_{\ov{\beta}}$ sur le disque $\Gamma'_{\ov{\beta}}$ de la manière suivante : à chaque arc de cercle obtenu à partir du demi-cercle $\partial \Gamma_{\ov{\beta}}$ en effectuant une homothétie $(x,y) \mapsto (\lambda x,y)$ (pour $\lambda \in \;]0,1]$), associons le cercle obtenu à partir du cercle $\partial' \Gamma_{\ov{\beta}}$ en effectuant une homothétie $(x,y) \mapsto (\lambda x,\lambda (y-\beta) + \beta)$ (ce qui donne un cercle de rayon $\frac{1}{10 \lambda b^2}$) : voir la figure~\ref{fig:gamma_gamma'}
\bfg
\begin{picture}(0,0)%
\includegraphics{gamma_gammap.pstex}%
\end{picture}%
\setlength{\unitlength}{2763sp}%
\begingroup\makeatletter\ifx\SetFigFont\undefined%
\gdef\SetFigFont#1#2#3#4#5{%
  \reset@font\fontsize{#1}{#2pt}%
  \fontfamily{#3}\fontseries{#4}\fontshape{#5}%
  \selectfont}%
\fi\endgroup%
\begin{picture}(6398,3868)(461,-3726)
\put(592,-3649){\makebox(0,0)[lb]{\smash{{\SetFigFont{8}{9.6}{\rmdefault}{\mddefault}{\updefault}{\color[rgb]{0,0,0}$\Gamma_{\overline{\beta}}$}%
}}}}
\put(4194,-3649){\makebox(0,0)[lb]{\smash{{\SetFigFont{8}{9.6}{\rmdefault}{\mddefault}{\updefault}{\color[rgb]{0,0,0}$\Gamma'_{\overline{\beta}}$}%
}}}}
\end{picture}%

\caption{L'homéomorphisme entre $\Gamma_{\ov{\beta}}$ et $\Gamma'_{\ov{\beta}}$}
\label{fig:gamma_gamma'}
\efg
On a ainsi construit un homéomorphisme de $C$ sur $Z_{>0}$, qui envoie chaque voisinage d'un point $c$ de $\partial \ov{C}$ dans $C$ sur un voisinage du point $g(c)$ dans $Z_{>0}$. Si on considère de plus l'homéomorphisme $i_A$ de $A$ sur $Z_0$, ces deux applications définissent une unique application $\eta'$ de $\ov{C} \cup_g A$ sur $Z$. Cette application est bijective et continue d'après la proposition~\ref{pro:recollement_fonction} puisqu'on peut prolonger l'homéomorphisme de $C$ sur $Z_{>0}$ par $i_A \circ g : \partial \ov{C} \ra Z_0$ en une application continue de $\ov{C}$ sur $Z$, compatible avec $g$. Par ailleurs l'application $\eta'$ est un homéomorphisme en restriction aux images de $C$ et $A$ dans $\ov{C} \cup_g A$, et la remarque sur l'image des voisinages assure que cette application est ouverte. Ainsi l'application $\eta'$ est un homéomorphisme de $\ov{C} \cup_g A$ sur $Z$. \ep

On peut voir sur la figure~\ref{fig:espacez} un schéma représentant l'espace $Z$ : le cône $\ov{C}$ se recolle sur l'espace des anneaux hawaïens $A$. Les pointillés suggèrent la façon dont le cône \og tourne \fg\ afin de réaliser cela.
\bfg
\includegraphics{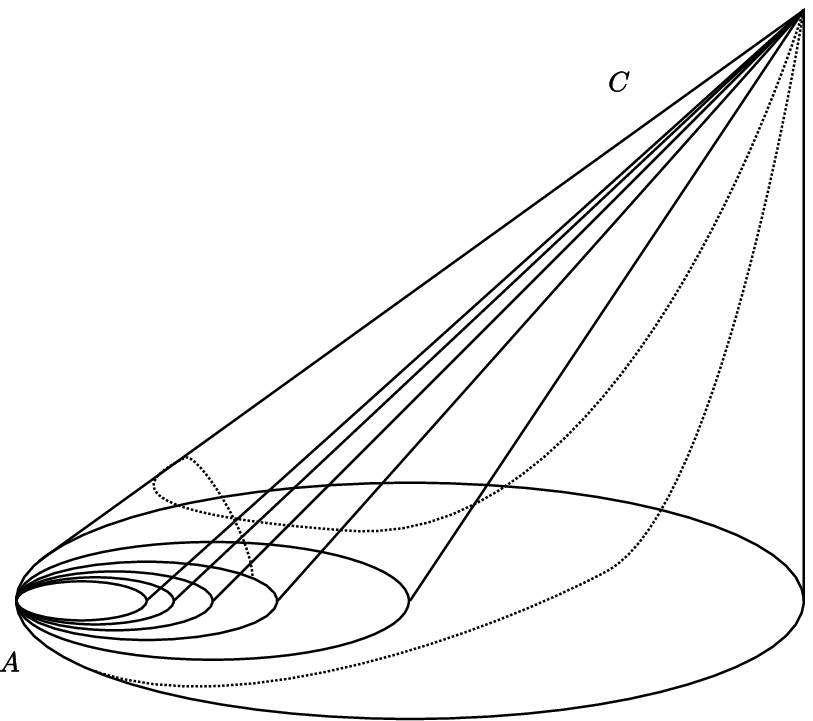}
\caption{L'espace $Z$}
\label{fig:espacez}
\efg
On a construit un homéomorphisme $\eta : Z \ra \ov{C} \cup_g A$. Pour $n \in \N \bs \{0\}$, si l'on compose cet homéomorphisme avec l'homéomorphisme
\beq A & \ra & A_n \\
\frac{1}{b}(1+e^{2i\theta}) & \mapsto & \frac{1}{nb}(1+e^{2i\theta}), \eeq
on obtient un plongement $\eta_n : Z \ra \ov{C} \cup_g A$, et son image est l'image de $C \cup A_n$ dans $\ov{C} \cup_g A$.

\bigskip

Considérons l'espace $C_\infty = \Ch(\Z) \times \ov{C}$, quotienté par la relation d'équivalence $(\{0\},\alpha,\ov{\beta}) \sim (\{0\},\alpha,\ov{\beta'})$ pour tous $\alpha \in [0,\infty]$ et $\ov{\beta},\ov{\beta'} \in \R/\Z$. Cet espace $C_\infty$ est constitué d'une suite de copies $\ov{C}_n$ du cône $\ov{C}$ s'accumulant sur leur axe $[0,\infty]$. Considérons l'application
\beq g_\infty : \partial C_\infty = \Ch(\Z) \times \partial \ov{C} & \ra & A \\
(n\Z,c) & \mapsto & \left\{ \begin{array}{ll} \frac{1}{nb}(1+e^{2i\theta}) & \mbox{ si } n \neq 0 \mbox{ et } g(c)=\frac{1}{b}(1+e^{2i\theta})\\ 0 & \mbox{ si } n = 0. \end{array} \right. \eeq

\bthm L'espace $\Ch(\R \times \Z)$ est homéomorphe au recollement $C_\infty \cup_{g_\infty} A$ des espaces $C_\infty$ et $A$ le long de $g_\infty$. \ethm

\bp Considérons l'application
\beq \phi':X = \Ch(\Z) \times Z & \ra & C_\infty \cup_{g_\infty} A \\
(n\Z,z) & \mapsto & \left\{ \begin{array}{ll} (n\Z,\eta_n(z)) \in C_\infty & \mbox{ si } n \neq 0\\
(\{0\},\alpha(z),0) \in C_\infty & \mbox{ si } n = 0. \end{array} \right.. \eeq

D'après la proposition précédente, l'application $\phi'$ est un plongement en restriction à $\{n\Z\} \times Z$, pour tout entier $n \in \N \bs \{0\}$. Puisque l'espace $\{n\Z\} \times Z$ est ouvert dans $X$ pour tout $n \in \N \bs \{0\}$, l'application $\phi'$ est un plongement en restriction à $(\Ch(\Z) \bs \{\{0\}\}) \times Z$. Par ailleurs il est clair que l'application $\phi'$ est continue en restriction à $\{\{0\}\} \times Z$. La continuité de l'application $\alpha$, et le fait que les cylindres $\ov{C_n}$ s'accumulent dans $C_\infty$, assurent que l'application $\phi'$ est continue sur $\{\{0\}\} \times Z_{>0}$. Par ailleurs le plongement $\eta_n$, pour $n \in \N \bs \{0\}$, assure que l'application $\phi'$ est continue sur $\{\{0\}\} \times Z_0$.

En conclusion, l'application $\phi'$ est une surjection continue de $X$ sur $C_\infty \cup_{g_\infty} A$. Il suffit de constater que la relation d'équivalence sur $X$ donnée par $\phi'(x) \sim \phi'(x')$ est exactement la même que celle considérée précédemment $\phi(x) \sim \phi(x')$, pour en conclure que les deux espaces quotients $C_\infty \cup_{g_\infty} A$ et $\Ch(\R \times \Z)$ sont homéomorphes. \ep

Nous avons donc démontré le théorème principal de cet article (et donc le théorème énoncé dans l'introduction). Voici une représentation de l'espace $\Ch(\R \times \Z)$ (figure~\ref{fig:espacex}).

\bfg
\includegraphics{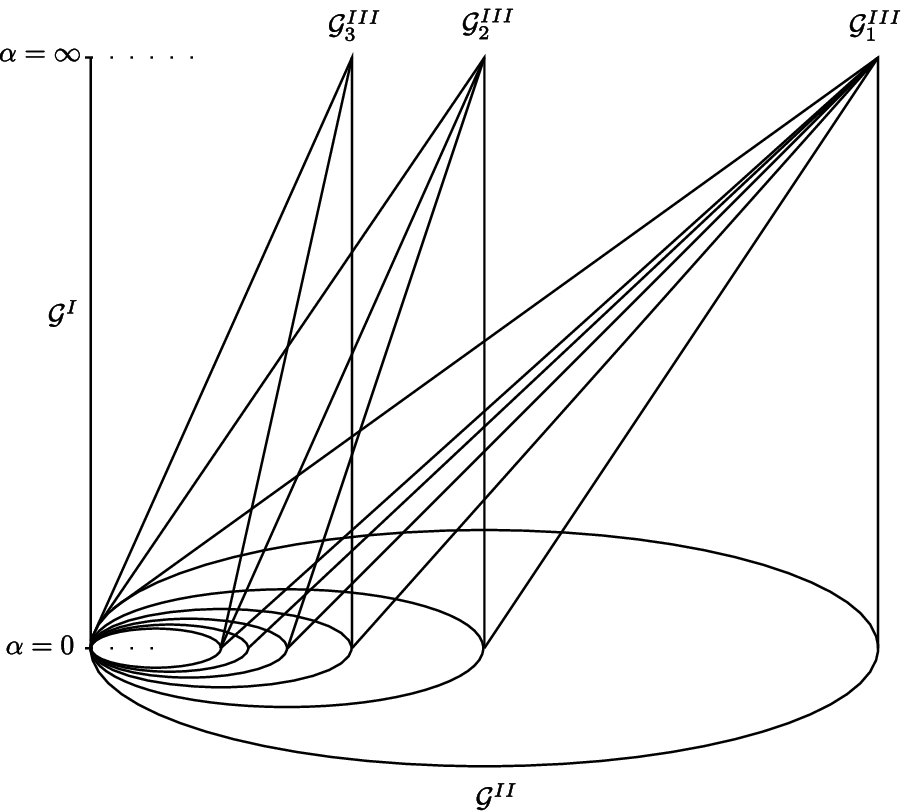}
\caption{L'espace $\Ch(\R \times \Z)$ des sous-groupes fermés de $\R \times \Z$}
\label{fig:espacex}
\efg

Remarquons, comme souligné par le rapporteur de cet article, le fait suivant.


\bpro Soit $\Ch_{\Z}(\R \times \Z)$ le sous-espace de $\Ch(\R \times \Z)$ constitué des sous-groupes cycliques infinis fermés de $\R \times \Z$, et $Q$ le quotient de $\R \times \Z \bs \{0,0\}$ par l'involution $(x,n) \mapsto (-x,-n)$. Alors l'application $\R \times \Z \bs \{0,0\}  \ra  \Ch_{\Z}(\R \times \Z)$ définie par $(x,n)  \mapsto \Z \cdot (x,n)$ induit un homéomorphisme $\eta$ de $Q$ sur $\Ch_{\Z}(\R \times \Z)$. \epro

\bp
Il est clair que l'application $\eta$ est bien définie, et bijective.




Soit $(x,n) \in \R \times \Z \bs \{0,0\}$, montrons que $\eta$ est continue en $\pm (x,n)$. Soit $(\pm (x_k,n_k))_{k \in \N}$ une suite de $Q$ convergeant vers $\pm (x,n)$. Choisissons les représentants tels que la suite $(x_k,n_k)_{k \in \N}$ converge vers $(x,n)$ dans $\R \times \Z$. Tout élément de $\Z \cdot (x,n)$ est clairement limite d'éléments de $(\Z \cdot (x_k,n_k))_{k \in \N}$. Réciproquement, soit $P$ une partie infinie de $\N$ et $(p_k)_{k \in P}$ une suite de $\Z$ telle que la suite $(p_kx_k,p_kn_k)_{k \in P}$ converge vers $(x',n')$. Puisque $(x,n) \neq (0,0)$, ceci impose que la suite $(p_k)_{k \in P}$ soit bornée : quitte à extraire, on peut donc supposer qu'elle est constante égale à $p \in \Z$. Ainsi $(x',n') = (px,pn) \in \Z \cdot (x,n)$. Ceci conclut la preuve de la continuité de l'application $\eta$.

Soit $\Z \cdot (x,n) \in \Ch_{\Z}(\R \times \Z)$, montrons que l'application $\eta^{-1}$ est continue en $\Z \cdot (x,n)$. Soit $(\Z \cdot (x_k,n_k))_{n \in \N}$ une suite de $\Ch_{\Z}(\R \times \Z)$ convergeant vers $\Z \cdot (x,n)$. Si la suite $(x_k,n_k)_{k \in \N}$ sortait de tout compact de $\R \times \Z$, la suite $(\Z \cdot (x_k,n_k))_{n \in \N}$ convergerait vers $\{0\}$. Si la suite $(x_k,n_k)_{k \in \N}$ convergeait vers $(0,0)$, la suite $(\Z \cdot (x_k,n_k))_{n \in \N}$ n'admettrait pas de limite discrète. Ainsi, quitte à extraire, on peut supposer que la suite $(x_k,n_k)_{k \in \N}$ converge vers $(x',n') \neq (0,0)$. Et d'après la première partie, dans ce cas nous avons $(x,n) = \pm (x',n')$. Donc l'application $\eta^{-1}$ est continue.
\ep

L'application $\eta$ a pour image $(\Ch^I \bs \{\R \times \{0\}\}) \cup \Ch^{II}$. Elle envoie la demi-droite $(\R \times \{0\}) / \pm$ sur la demi-droite $[0,\infty[\ \subset [0,\infty] \simeq \Ch^I$, par la valeur absolue. Et elle envoie la droite $\R \times \{\pm n\}$, pour $n \in \N \bs \{0\}$, sur le cercle épointé $A_n \bs \{0\} \subset A \simeq \Ch^{II}$, par l'application $(x,\pm n) \mapsto \frac{1}{n}(1+e^{2i\tan x})$. Ainsi, l'espace $\Ch_{\Z}(\R \times \Z)$ est homéomorphe au recollement d'un intervalle $[0,\infty[$ sur les anneaux hawaïens $A$, en identifiant $0 \in [0,\infty[$ et $0 \in A$.

\subsection{Le groupe fondamental de l'espace $\Ch(\R \times \Z)$}

Nous allons maintenant nous intéresser au groupe fondamental de l'espace $\Ch(\R \times \Z)$.

Considérons les sous-espaces $Z_0=\{\alpha = 0\}$, $Z_{>0}=\{\alpha > 0\}$ et $Z_{\leq 1}=\{\alpha \in [0,1]\}$ de $Z$. Puis les espaces $\widetilde{X_0}$, $\widetilde{X_{>0}}$ et $\widetilde{X_{\leq 1}}$, images de $\Ch(\Z) \times Z_0$, $\Ch(\Z) \times Z_{>0}$ et $\Ch(\Z) \times Z_{\leq 1}$ dans $\widetilde{X}$.

\bpro L'espace $\widetilde{X_{\leq 1}}$ se rétracte par déformation forte sur $\widetilde{X_0}$. \epro

\bp Considérons l'application $h: \widetilde{X_{\leq 1}} \times [0,1] \ra \widetilde{X_{\leq 1}}$, induite par passage au quotient de
\beq \Ch(\Z) \times Z_{\leq 1} \times [0,1] & \ra & \Ch(\Z) \times Z_{\leq 1} \\
(n\Z,z,t) & \mapsto & \left\{ \begin{array}{ll} (n\Z,\eta_n^{-1}(t\alpha,\ov{\beta})) & \mbox{ si } n \neq 0 \mbox{ et } \eta_n(z) = (\alpha,\ov{\beta}) \\ (\{0\},t\alpha) & \mbox{ si } n = 0 \mbox{ et } \alpha=\alpha(z). \end{array} \right. \eeq
Cette application est continue, fixe $\widetilde{X_0}$ pour tout $t \in [0,1]$, vaut l'identité pour $t=1$ et a pour image $\widetilde{X_0}$ pour $t=0$. Par conséquent, l'espace $\widetilde{X_{\leq 1}}$ se rétracte par déformation forte sur $\widetilde{X_0}$. \ep

Le résultat suivant est bien connu, voir par exemple~\cite[Example~1.25, p.~49]{hatcher}.

\bpro Le premier groupe d'homologie entière $H_1(A)$ de l'espace $A$ n'est pas dénombrable. \hfill $\square$ \epro

\bthm Le premier groupe d'homologie entière, et donc également le groupe fondamental, de l'espace $\Ch(\R \times \Z)$ est non dénombrable. \ethm

\bp
L'espace $\widetilde{X_{\leq 1}}$ se rétracte par déformation forte sur $\widetilde{X_0}$, espace qui est homéomorphe à $A$. Ainsi le groupe $H_1(\widetilde{X_{\leq 1}})$ n'est pas dénombrable.

L'espace $\widetilde{X_{>0}}$ est homéomorphe à une suite de disques ouverts disjoints $D_n$ s'accumulant sur un rayon $D_0$. D'après la proposition~\cite[Proposition~2.6, p.~109]{hatcher}, nous avons $H_1(\widetilde{X_{>0}}) = 0$.

Enfin remarquons que l'espace $\widetilde{X_{>0}} \cap \widetilde{X_{\leq 1}}$ est homéomorphe à une suite de cylindres disjoints $C_n$ s'accumulant sur leur axe $C_0$. D'après la même proposition, nous avons $H_i(\widetilde{X_{>0}} \cap \widetilde{X_{\leq 1}}) \simeq \oplus_{n \in \N} H_i(C_n)$ pour tout entier $i \in \N$. Ainsi
\beq H_0(\widetilde{X_{>0}} \cap \widetilde{X_{\leq 1}}) & \simeq & \bigoplus_{n \in \N} \Z \\
\mbox{et } H_1(\widetilde{X_{>0}} \cap \widetilde{X_{\leq 1}}) & \simeq & \bigoplus_{n \in \N \bs \{0\}} \Z .\eeq
En particulier, ces groupes sont dénombrables.

\bigskip

Puisque les intérieurs des espaces $\widetilde{X_{>0}}$ et $\widetilde{X_{\leq 1}}$ recouvrent $\widetilde{X}$, on peut écrire la suite exacte de Mayer-Vietoris pour l'homologie entière (voir~\cite[pp.~149--153]{hatcher}) :
$$ \ldots H_2(\widetilde{X}) \ra H_1(\widetilde{X_{>0}} \cap \widetilde{X_{\leq 1}}) \ra H_1(\widetilde{X_{>0}}) \oplus H_1(\widetilde{X_{\leq 1}}) \stackrel{\gamma}{\ra} H_1(\widetilde{X}) \ra H_0(\widetilde{X_{>0}} \cap \widetilde{X_{\leq 1}}) \ldots $$

En conclusion, le morphisme $\gamma : H_1(\widetilde{X_{\leq 1}}) \ra H_1(\widetilde{X})$ a un noyau $H_1(\widetilde{X_{>0}} \cap \widetilde{X_{\leq 1}})$ et un conoyau $H_0(\widetilde{X_{>0}} \cap \widetilde{X_{\leq 1}})$ dénombrables. Or l'espace de départ $H_1(\widetilde{X_{\leq 1}})$ est non dénombrable, donc l'espace d'arrivée $H_1(\widetilde{X}) \simeq H_1(\Ch(\R \times \Z))$ n'est pas dénombrable. \ep

\bigskip

\sign

\bibliographystyle{perso}
\bibliography{biblio}

\end{document}